\documentclass[12pt,oneside]{amsproc}

\usepackage[cp1251]{inputenc}
\usepackage[russian]{babel}
\usepackage{amssymb,amsmath}
\numberwithin{equation}{section}

\newtheorem{tm}{Теорема}[section]

\newtheorem{stat}{Утверждение}[section]
\newtheorem{rem}{Замечание}[section]

\newcommand{\Wo}{{\raisebox{0.2ex}{\(\stackrel{\circ}{W}\)}}{}}

\title[]{О спектре оператора Якоби с экспоненциально растущими матричными элементами}
\author{
И.~А.~Шейпак\footnote{%
Работа выполнена при поддержке грантов РФФИ \No\,07-01-00283,  поддержки ведущих научных школ \No\,НШ-2372.2008.1 и российско-украинского гранта РФФИ \No\,09-01-90408}%
\address{Московский государственный университет
им.~М.~В.~Ломоносова, механико-математический факультет}
\email{iasheip@mech.math.msu.su}}

\begin{document}
\noindent УДК~517.984
\begin{abstract}
В статье рассматривается класс матриц Якоби с быстро растущими матричными элементами. В пространстве квадратично суммируемых с некоторым весом
последовательностей этой матрице отвечает симметрический оператор. Доказывается, что задача на собственные значения некоторого самосопряжённого
расширения этого оператора эквивалентна задаче на собственные значения оператора Штурма--Лиувилля с дискретным самоподобным самоподобным весом.
Находятся асимптотические формулы для собственных значений.
\end{abstract}
\maketitle

\section{Введение}
Исследуется задача на собственные значения трёхдиагональной якобиевой матрицы вида

\begin{equation}\label{matrix}
\begin{pmatrix}
\alpha & \beta & 0 & 0 & \dots & 0 & \dots\\
\gamma & \alpha q & \beta q & 0 & \dots & 0 & \dots\\
 0 &  \gamma q & \alpha q^2 & \beta q^2 & \dots & 0 & \dots\\
\dots\\
0 & 0 & 0 & \gamma q^{n-1} & \alpha q^n & \beta q^n & 0 \\
\dots & \dots & \dots & \dots & \dots & \dots & \dots
\end{pmatrix},
\end{equation}
где $q>1$.

В работе \cite{TUR} были изучены спектральные свойства класса двухдиагональных симметричных матриц (на главной диагонали стоят нули), в которых
внедиагональные элементы $b_n$ имеют вид $q^{n^s}$, $0<q<1$, $s$ --- произвольное натуральное число. Там же было доказано, что
$b_{2n+1}<\lambda_n<b_{2n-1}$. При этом использовалась достаточно тяжёлая техника, основанная на свойствах целых трансцендентных функций.

В работе~\cite{KOZHAN} был изучен более широкий класс двухдиагональных симметричных матриц Якоби, в которых внедиагональные элементы $b_n$
положительны, и для них выполнено условие
$$
\lim_{n\to\infty}\dfrac{b_{n+1}}{b_n}=0.
$$
В этом случае получены не такие точные оценки на собственные значения, как в работе \cite{TUR}. Асимптотические формулы имеют вид
$\lambda_n=b_{2n-1}(1+o(1))$,  $n\to\infty$. Метод получения асимптотических формул опирался на исследование квадратичной формы соответствующего
оператора.

Заметим, что в нашем случае, $q>1$ и оператор с такой матрицей неограничен. Необходимо дополнительно описывать его область определения. Также
важен вопрос об индексах дефекта такого оператора. В случае индекса дефекта $(1,1)$ встаёт вопрос о самосопряжённых расширениях такого
оператора.

В дальнейшем мы всегда будем считать, что $\alpha\beta>0$. Это условие достаточно часто выполняется для матриц, возникающих в различных
приложениях. Для матриц конечного порядка это условие приводит к понятию \emph{нормальной матрицы} (см., например, \cite{GANTKR}).

Для данной работы метод, на основе которого выводятся асимптотические формулы для собственных значений матрицы~\eqref{matrix}, заключается в
том, что задача на собственные значения для оператора~\eqref{matrix} оказывается эквивалентной задаче на собственные значения задачи
Штурма--Лиувилля с самоподобным дискретным весом. В связи с этим, в \S 2 мы приведём некоторые сведения о самоподобных функциях, порождающих
дискретные веса (так называемые \emph{самоподобные функции нулевого спектрального порядка}, подробнее см. \cite{SH2}, \cite{VLASH3}). В \S 3
задача Штурма--Лиувилля с дискретным самоподобным весом сводится к матричному виду.

\section{Самоподобные функции нулевого спектрального порядка}
Пусть числа $a\in(0,1)$ и $d$ удовлетворяют условию
\begin{equation}\label{szhim}
a|d|^2<1.
\end{equation}
Определим в $L_2[0,1]$ оператор  $G$, действующий по правилу
\begin{equation}\label{operator podobia}
G(f)(x)=\beta_1\cdot\chi_{[0,1-a)}(x)+\left(d\cdot f\left(\frac{x-1+a}{a}\right)+\beta_2\right)\cdot\chi_{(1-a,1]}(x),
\end{equation}
где $\beta_1$, $\beta_2$ --- произвольные действительные числа.

Условие~\eqref{szhim} влечёт, что оператор $G$ является сжимающим в $L_2[0,1]$ (см. лемму 3.1,~\cite{VLASH1}) и, следовательно, имеет
единственную неподвижную точку. Функция, удовлетворяющая уравнению $G(P)=P$, где отображение $G$ задаётся соотношением~\eqref{operator podobia},
является самоподобной функцией нулевого спектрального порядка. Набор чисел $a$, $d$, $\beta_1$ и $\beta_2$ называется параметрами самоподобия,
задающими функцию $P$.

Из определения следует, что функция $P$ является кусочно постоянной, причем принимает значения
\begin{gather}\label{samopod_function1}
P(x)=\beta_1, \text{ при } x\in[0,1-a),\\
\label{samopod_function2}
P(x)=d^{k}\beta_1+\beta_2(1+d+\ldots+d^{k-1}), \text{ при } x\in(1-a^k,1-a^{k+1}),\quad k=1,2,\ldots.
\end{gather}

Более подробно о самоподобных функциях в пространствах $L_p[0,1]$ см.\cite{SH1}. Общая конструкция и свойства самоподобных функции нулевого
спектрального порядка изучены в~\cite{SH2}, см. также~\cite{VLASH3}.

\section{Задача Штурма-Лиувилля в сингулярным самоподобным весом}
Рассмотрим следующую граничную задачу
\begin{gather}\label{sturm_weight_problem}
-y''-\lambda\rho y=0,\\\label{sturm_weight_uslovia} y(0)=y(1)=0,
\end{gather}
где $\rho=P'$ (c смысле обобщённых функций), т.е.
\begin{equation}\label{weight}
\rho:=\sum_{k=1}^\infty m_k\delta(x-(1-a^k)),
\end{equation}
где $m_k=d^{k-1}(d\beta_1+\beta_2-\beta_1)$. Из условия $P\in L_2[0,1]$ следует, что $\rho\in\Wo^{-1}_2[0,1]$.

В дальнейшем через $\mathfrak H$ мы будем обозначать пространство $\Wo^1_2[0,1]$, снабжённое скалярным произведением
$$
\langle y,z\rangle =\int\limits_0^1y'\overline{z'}dx
$$
Через $\mathfrak H'$ мы будем обозначать пространство, двойственное к $\mathfrak H$ относительно $L_2[0,1]$, т.е. получаемое пополнением
пространства $L_2[0,1]$ по норме
$$
\|y\|_{\mathfrak H'}=\sup\limits_{\|x\|_{\mathfrak H'}}\left|\int\limits_0^1y\overline{z}dx\right|.
$$

Как и в работе~\cite{VLASH3} рассмотрим оператор вложения $J:\mathfrak H\to L_2[0,1]$. Непосредственно из определения пространства
$\mathfrak H'$ вытекает возможность непрерывного продолжения сопряжённого оператора $J^*:L_2[0,1]\to \mathfrak H$ до изометрии
$J^+:\mathfrak H'\to \mathfrak H$.

Как и в предшествующих работах~\cite{VLASH1}, \cite{VLASH2} и \cite{VLASH3} в качестве операторной модели
задачи~\eqref{sturm_weight_problem}--\eqref{sturm_weight_uslovia} мы будем рассматривать линейный пучок $T_\rho:\mathfrak H\to\mathfrak H'$
ограниченных операторов, удовлетворяющий тождеству
\begin{equation}\label{eq:difur_pichok}
\forall \lambda, \forall y\in \mathfrak H \qquad \langle
J^+T_\rho(\lambda)y,y\rangle =\int\limits_0^1\left(|y'|^2+\lambda P\cdot(|y|^2)'\right)dx.
\end{equation}

Известно (теорема 4.1,~\cite{VLASH1}), что спектр задачи чисто дискретен при любом весе из пространства $\Wo^{-1}_2[0,1]$.

Рассмотрим задачу~\eqref{sturm_weight_problem}--\eqref{sturm_weight_uslovia} при условии, что вес является обобщённой производной самоподобной
функции $P$ нулевого спектрального порядка, т.е. когда вес определяется соотношением~\eqref{weight}.

\subsection{Собственные функции задачи Штурма--Лиувилля с весом, являющимся обобщённой производной функции нулевого спектрального порядка}

Т.к. функция $P$ кусочно-постоянная,  собственные функции этой задачи можно искать в виде кусочно-линейной функции $y\in \Wo_2[0,1]$, заданной
на каждом промежутке $(1-a^{k},1-a^{k+1})$ $k=0,1,2,\ldots$ формулами
\begin{gather}\label{eigenfunctions1}
y(x)=s_1 x, \quad x\in[0,1-a),\\
\label{eigenfunctions2}
y(x)=s_k x+t_k, \quad x\in (1-a^{k-1},1-a^{k}), \quad k=2,3,\ldots.
\end{gather}

Условия непрерывности функции $y$ в точках $x_k:=1-a^k$, $k=1,2,\ldots$ означают что выполнены уравнения
\begin{equation}\label{nepreryvnost}
s_{k-1}(1-a^{k-1})+t_{k-1}=s_{k}(1-a^{k-1})+t_{k}, \quad k=2,3,\ldots.
\end{equation}

Подстановка формул~\eqref{eigenfunctions1}--\eqref{eigenfunctions2} в уравнение~\eqref{sturm_weight_problem} с учётом
уравнений~\eqref{nepreryvnost} приводит к следующей системе уравнений на числа
$s_k$:
\begin{gather*}
s_1-s_2=\lambda(d\beta_1+\beta_2-\beta_1)(1-a)s_1,\\
s_2-s_3=\lambda\left(d^2\beta_1+d(\beta_2-\beta_1)\right)(1-a)(s_1+as_2),\\
s_3-s_3=\lambda\left(d^3\beta_1+d^2(\beta_2-\beta_1)\right)(1-a)(s_1+as_2+a^2s_3),\\
\ldots,\\
s_n-s_{n+1}=\lambda\left(d^n\beta_1+d^{n-1}(\beta_2-\beta_1)\right)(1-a)(s_1+as_2+\ldots+a^{n-1}s_n),\\
\ldots
\end{gather*}

Введём следующие обозначения:
$$
r=(1-a)(d\beta_1+\beta_2-\beta_1), \quad q=\dfrac{1}{ad}.
$$

Из условия~\eqref{szhim} следует, что
$$
|q|>1.
$$

В дальнейшем будем считать, что параметры самоподобия таковы, что $r\ne 0$, а $d>0$. Следовательно, $q>1$.

Краевое условие при $x=0$ выполнено в силу формулы~\eqref{eigenfunctions1}. Необходимо учесть краевое условие при $x=1$. Из
формулы~\eqref{eigenfunctions2} следует, что
$$
y(1)=\lim_{n\to\infty}(s_n+t_n).
$$

Из~\eqref{nepreryvnost} несложно получить, что
$$
t_n=\sum_{k=1}^{n-1}s_k(x_k-x_{k-1})-s_nx_{n-1}.
$$

Окончательно получаем, что
$$
y(1)=(1-a)\sum_{k=1}^\infty a^{k-1}s_k=0,
$$
т.е. последовательность $\{s_k\}_{k=1}^\infty$ удовлетворяет условию
\begin{equation}\label{kraev_usl_v_1}
\sum_{k=1}^\infty a^{k-1}s_k=0.
\end{equation}

Заметим, что $\|y\|_{\mathfrak H}^2=\sum\limits_{k=1}^\infty (1-a^{k}-(1-a^{k-1}))s_k^2=(1-a)\sum\limits_{k=1}^\infty a^{k-1}s_k^2$.  Таким
образом, задача~\eqref{sturm_weight_problem}--\eqref{sturm_weight_uslovia} с дискретным самоподобным весом эквивалентна следующей задаче в
пространстве последовательностей $\{s_k\}_{k=1}^\infty$, суммируемых в квадрате с весом $\omega=\{\omega_k\}_{k=1}^\infty$, где
$\omega_k=a^{k-1}$, и удовлетворяющих условию~\eqref{kraev_usl_v_1}:

\begin{multline}\label{puchok}
\begin{pmatrix}
1 & -1 &  0 &  0 & \dots & 0 & \dots\\
0 &  1 & -1 &  0 & \dots & 0 & \dots\\
0 &  0 &  1 & -1 & \dots & 0 & \dots\\
\dots\\
0 & 0 & 0 & 0 & 1 & -1 & 0 \\
\ldots & \ldots & \ldots & \ldots & \ldots & \ldots & \ldots
\end{pmatrix}
\begin{pmatrix}
s_1\\
s_2\\
s_3\\
\ldots\\
s_k\\
\ldots
\end{pmatrix}=\\=
\lambda r
\begin{pmatrix}
1   &  0    &  0 &  0 & \dots & 0     & \dots\\
d   &  da   &  0 &  0 & \dots & 0     & \dots\\
d^2 &  d^2a &  (da)^2 & 0     & \dots & 0 & \dots\\
\dots\\
d^{k-1} & d^{k-1}a & d^{k-1}a^2 & \ldots & (da)^{k-1} & 0 & \ldots \\
\ldots & \ldots & \ldots & \ldots & \ldots & \ldots & \ldots
\end{pmatrix}\begin{pmatrix}
s_1\\
s_2\\
s_3\\
\ldots\\
s_k\\
\ldots
\end{pmatrix}
\end{multline}

Введём в рассмотрение следующие операторы, заданные матрицами:
$$
A=\begin{pmatrix}
1 & -1 &  0 &  0 & \dots & 0 & \dots\\
0 &  1 & -1 &  0 & \dots & 0 & \dots\\
0 &  0 &  1 & -1 & \dots & 0 & \dots\\
\dots\\
0 & 0 & 0 & 0 & 1 & -1 & 0 \\
\ldots & \ldots & \ldots & \ldots & \ldots & \ldots & \ldots
\end{pmatrix}\quad
B=\begin{pmatrix}
1   &  0    &  0 &  0 & \dots & 0     & \dots\\
d   &  da   &  0 &  0 & \dots & 0     & \dots\\
d^2 &  d^2a &  (da)^2 & 0     & \dots & 0 & \dots\\
\dots\\
d^{k-1} & d^{k-1}a & d^{k-1}a^2 & \ldots & (da)^{k-1} & 0 & \ldots \\
\ldots & \ldots & \ldots & \ldots & \ldots & \ldots & \ldots\\
\end{pmatrix}.
$$
Таким образом задача \eqref{eq:difur_pichok} эквивалентна задаче
$$
As=\lambda rBs,
$$
рассматриваемой в пространстве последовательностей $l_{2,a}$, удовлетворяющих условию \eqref{kraev_usl_v_1}.

Несложно проверить, что обратный к оператору $B$  имеет вид
$$
B^{-1}=\begin{pmatrix}
1 & 0 & 0 & 0 & \ldots\\
-dq & q & 0 & 0 & \ldots\\
0 & -dq^2 & q^2 & 0 & \ldots\\
0 & 0 & -dq^3 & q^3 & \ldots\\
\ldots & \ldots & \ldots & \ldots & \ldots \\
\end{pmatrix}.
$$
Возникает естественный вопрос, является ли спектральная задача~\eqref{puchok} эквивалентна одной из следующих задач
$$
B^{-1}As=\lambda r s\qquad \text{ или } \qquad AB^{-1} u=\lambda r u\quad (\text{здесь }u:=Bs)
$$
в каком-нибудь пространстве последовательностей?

Рассмотрим вещественное число $w\ne 0$. Обозначим пространство последовательностей $\{v_k\}_{k=1}^\infty$, удовлетворяющих условию
$$
\sum_{k=1}^\infty w^{k-1}v_k^2<\infty,
$$
через $l_{2,w}$. Скалярное произведение в этом пространстве будем обозначать через $\langle \cdot,\cdot \rangle_{w}$.

\begin{stat}
Оператор $AB^{-1}$ симметричен в $l_{2,1/d}$. 
\end{stat}
\begin{proof}
Непосредственным вычислением несложно убедится, что
$$
\langle AB^{-1} u,v \rangle_{1/d}=\langle u, AB^{-1} v\rangle_{1/d}=(1+dq)\sum_{k=1}^\infty\left(\dfrac{q}{d}\right)^{k-1} u_k v_k
-q\sum_{k=1}^\infty \left(\dfrac{q}{d}\right)^{k-1}(u_{k+1}v_k+u_kv_{k+1}).
$$
\end{proof}

\begin{stat}
Область определения сопряжённого оператора к $AB^{-1}$ состоит из всех последовательностей $u\in l_{2,1/d}$, таких что
\begin{equation}\label{eq:domain_sopr}
\sum_{k=2}^\infty \dfrac{1}{d^{k-1}}\left(-dq^{k-1}u_{k-1}+(1+dq)q^{k-1}u_k+q^ku_{k+1}\right)^2<\infty.
\end{equation}
\end{stat}
\begin{proof}
Заметим, что матрица $AB^{-1}$ имеет вид
\begin{equation}\label{eq:matr_jak}
AB^{-1}=\begin{pmatrix} 1+dq & -q & 0 & 0 & \ldots\\
-dq & (1+dq)q & -q^2 & 0 & \ldots\\
0 & -dq^2 & (1+dq)q^2 & -q^3 & \ldots\\
\ldots  & \ldots  & \ldots  & \ldots  & \ldots
\end{pmatrix}.
\end{equation}
Для завершения доказательства утверждения осталось применить \cite{AHIEZER} (п. 1.1, стр. 174).
\end{proof}

Можно определить индексы дефекта оператора $AB^{-1}$. Задача \eqref{sturm_weight_problem}--\eqref{sturm_weight_uslovia} самосопряжена. Условие
\eqref{eigenfunctions1} учитывает только одно краевое условие (в нуле). Остаётся ещё одно условие на собственные функции в единице.
Следовательно, индексы дефекта оператора $AB^{-1}$ равны $(1,1)$. Несложно проверить, что условие \eqref{kraev_usl_v_1} на последовательность
$s\in l_{2,a}$ переходит в условие
\begin{equation}\label{eq:kraev_usl_u}
\lim_{n\to\infty}\dfrac{u_n}{d^{n-1}}=0
\end{equation}
на последовательность $u:=Bs$ в пространстве $l_{2,1/d}$.

\begin{tm}\label{tm:equiv}
Задача \eqref{eq:difur_pichok}, эквивалента задаче
$$
AB^{-1}u=\lambda r u
$$
в пространстве $l_{2,1/d}$ с условием \eqref{eq:domain_sopr}, \eqref{eq:kraev_usl_u} на последовательность $u=(u_1,u_2,\ldots)$, где $u=Bs$,
$s\in l_{2,a}$.
\end{tm}
\begin{proof}
Представление \eqref{weight} веса $\rho$, а также условие \eqref{eigenfunctions1} приводит к тому, что квадратичная форма задачи
\eqref{eq:difur_pichok} в пространстве $l_{2,a}$ принимает вид
\begin{equation}\label{eq:diskr_kvadr_forma}
\sum_{k=1}^\infty a^{k-1}s_k=\lambda (1-a)\beta \sum_{k=1}^\infty d^{k-1}\left(\sum_{j=1}^k a^{j-1}s_j \right)^2.
\end{equation}
Сведём задачу $AB^{-1} u=\lambda r  u$ к её квадратичной форме $\langle AB^{-1} u,u\rangle_{1/d}=\lambda r \langle u, u\rangle_{1/d}$, которую
можно также переписать в виде
$$
\langle As,Bs\rangle_{1/d}=\lambda r \langle Bs, Bs\rangle_{1/d}.
$$
Используя вид операторов $A$ и $B$ убеждаемся, что эта квадратичная форма также имеет вид \eqref{eq:diskr_kvadr_forma}.
\end{proof}

\begin{rem}
Таким образом, одним из самосопряжённых расширений симметрического оператора $AB^{-1}$ в пространстве $l_{2,1/d}$ является оператор, заданный
той же матрицей, а его области определения принадлежат те последовательности $u\in l_{2,1/d}$, которые удовлетворяют условиям
\eqref{eq:domain_sopr} и \eqref{eq:kraev_usl_u}. Задача на собственные значения этого самосопряжённого расширения эквивалентно задаче
\eqref{sturm_weight_problem}--\eqref{sturm_weight_uslovia} (или задаче \eqref{eq:difur_pichok}).
\end{rem}

В следующей теореме мы подразумеваем именно это самосопряжённое расширение и обозначаем его через  $L$. Заметим также, что якобиева матрица
$AB^{-1}$ принадлежит классу матриц вида~\eqref{matrix} ($\alpha=1+dq=1+\dfrac{1}{a}$, $\beta=-q$, $\gamma=-dq=-\dfrac{1}{a}$) (см.
\eqref{eq:matr_jak}).

\begin{tm}
Существует такое положительное число $с$, что для собственных значений оператора $L$, занумерованных в порядке возрастания, справедлива
асимптотическая формула при $k\to\infty$

\begin{equation}\label{asymp}
\lambda_k=c q^k(1+o(1)).
\end{equation}
\end{tm}
\begin{proof}
Утверждение теоремы следует из теоремы \ref{tm:equiv} и работы (\cite{VLASH3}, Теорема 4.1).
\end{proof}

\subsection{Индефинитный случай}
Можно рассматривать задачу \eqref{sturm_weight_problem}--\eqref{sturm_weight_uslovia} и в случае, если число $d<0$. В этом случае оператор $L$ с
областью определения \eqref{eq:domain_sopr}, \eqref{eq:kraev_usl_u} будет  самосопряжённым в пространстве с индефинитной метрикой.

А именно, определим операторы ортогонального проектирования в пространстве $l_{2,1/d}$: $P_+: e_k\to e_{k}$, $k=1,3,\ldots, 2n-1,\ldots$,
$P_+:e_k\to 0$, $k=2,4,\ldots, 2n,\ldots$; $P_-: e_k\to 0$, $k=1,3,\ldots, 2n-1,\ldots$, $P_-:e_k\to e_k$, $k=2,4,\ldots,
2n,\ldots$$n\in\mathbb{N}$. Определим также оператор $J=P_+-P_-$. Несложно проверить, что оператор $L$ будет самосопряжённым в $J$-метрике.
Кроме того, из теоремы 4.3( \cite{VLASH3}) следует следующее утверждение
\begin{tm}
Пусть $d<0$ и $r\ne 0$. Тогда существует такое число $c>0$, что для положительных собственных значений $\{\lambda_k\}_{k=1}^\infty$ оператора
$L$, занумерованных в порядке возрастания, справедлива асимптотическая формула
$$
\lambda_{k+1}=c q^{2k}(1+o(1)),
$$
а для отрицательных собственных значений $\{\lambda_{-k}\}_{k=1}^\infty$ оператора $L$, занумерованных в порядке возрастания, справедлива
асимптотическая формула
$$
\lambda_{-(k+2)}=-c q^{2k+1}(1+o(1)).
$$
\end{tm}

\end{document}